# IGNORABILITY FOR CATEGORICAL DATA

By Manfred Jaeger

*Aalborg Universitet*

We study the problem of ignorability in likelihood-based inference from incomplete categorical data. Two versions of the coarsened at random assumption (*car*) are distinguished, their compatibility with the parameter distinctness assumption is investigated and several conditions for ignorability that do not require an extra parameter distinctness assumption are established.

It is shown that *car* assumptions have quite different implications depending on whether the underlying complete-data model is saturated or parametric. In the latter case, *car* assumptions can become inconsistent with observed data.

**1. Introduction.** In a sequence of papers Rubin [15], Heitjan and Rubin [11] and Heitjan [9, 10] have investigated the question under what conditions a mechanism that causes observed data to be incomplete or, more generally, *coarse*, can be ignored in the statistical analysis of the data. The key condition that has been identified is that the data should be *missing at random* (mar), respectively, *coarsened at random* (*car*). Similar conditions were independently proposed by Dawid and Dickey [4]. A second condition needed in Rubin's [15] derivation of ignorability is *parameter distinctness* (*pd*).

A case of particular practical interest is the one of incomplete or coarse categorical data. Traditionally associated with the analysis of contingency tables in terms of log-linear models, categorical data today also plays an important role in learning probabilistic models for artificial intelligence applications [12]. For these applications graphical models or Bayesian networks are used [2, 3, 13]. Incomplete data here is particularly prevalent, and the analysis of Rubin and Heitjan is widely cited in the field.

In this paper we take a closer look at the way ignorability is established for likelihood-based inference through the *car* and *pd* assumptions. It is found that one has to distinguish a weak version of *car* that is given as a condition on the joint distribution of complete and coarse data, and a strong version of *car* that is given as a condition on the conditional distribution of the coarse data. The two versions of *car* lead to quite different theoretical results and practical implications for likelihood-based inference. We consider in detail the dependencies between the *car* and the *pd* assumptions, and find that for weak *car* these two assumptions are incompatible unless further assumptions on the parameter of interest, or on the coarsening process, are made. In contrast, *pd* is implied by strong *car* (Section 3). For the case of an underlying saturated complete-data model ignorability results can be derived from weak *car* alone without making the *pd* assumption. Our main result identifies the maxima of the observed-data likelihood under either *car* assumption as exactly those complete-data distributions that are compatible with the *car* assumption and the observed data (Section 4.1). For nonsaturated complete-data models no analogous results hold. Even for very simple parametric models *car* becomes a testable assumption that can be rejected against an alternative hypothesis (Section 4.2).

**2. Coarse data models.** We use a very general and abstract model for categorical data: complete data is taken to consist of realizations $x_1, \ldots, x_N$ of independent identically distributed random variables $X_1, \ldots, X_N$ that take values in a finite set $W = \{w_1, \ldots, w_n\}$. The $w_i$ can be the cells of a contingency table, for instance. The distribution of the $X_i$ is assumed to belong to a parametric family $\{P_\theta | \theta \in \Theta\}$, where $\Theta \subseteq \mathbb{R}^k$ for some $k \in \mathbb{N}$. For this paper the analytic form of a parametric family will not be important, and only the subset of distributions contained in the family is relevant. For that reason we may generally assume that

$$\Theta \subseteq \Delta^n := \Big\{(p_1, \ldots, p_n) \in [0,1]^n \Big| \sum p_i = 1\Big\}$$

with

$$P_\theta(w_i) = p_i, \qquad \theta = (p_1, \ldots, p_n) \in \Theta.$$

Any $\Theta \subseteq \Delta^n$ is called a *complete-data* model. $\Theta = \Delta^n$ is the *saturated* complete-data model. In the saturated model, as well as in most of the important parametric models for categorical data (e.g., log-linear models), different parameters $\theta, \theta'$ may define distributions $P_\theta, P_{\theta'}$ with different sets of support. Most of the results of this paper address difficulties that arise out of this.

When data is incomplete, then the exact value $x_i$ of $X_i$ is not observed. According to the general coarse data model of Heitjan and Rubin [11] one observes instead a subset $U_i$ of $W$. More specifically, Heitjan and Rubin model coarse data by introducing additional coarsening variables $G_i$, and taking $U_i$ to be a function $Y(x_i, g_i)$ of the complete data $x_i$ and the value $g_i$ of the coarsening variable. In the following definition we take a slightly different approach, and model the coarsening process directly by a joint distribution of $X_i$ and the observed coarse data $U_i$. For categorical data this is simpler, and avoids a sometimes artificial construction of a suitable coarsening variable.

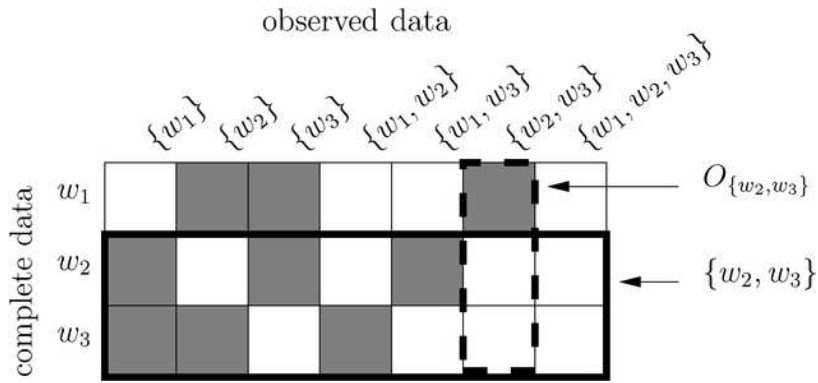

Fig. 1. *Coarse data space.*

DEFINITION 2.1. Let $W = \{w_1, \ldots, w_n\}$. The *coarse data space* for $W$ is
$$\Omega(W) := \{(w, U) | w \in W, U \subseteq W : w \in U\}.$$

When specific reference to $W$ is not needed, we write $\Omega$ for $\Omega(W)$. An element $(w, U) \in \Omega$ stands for the event that the true value of a $W$-valued random variable $X$ is $w$ and the coarse value $U$ is observed. A subset $U \subseteq W$ defines two different subsets in $\Omega$: $O_U := \{(w, U) \in \Omega | w : w \in U\}$, which is the event that $U$ is observed, and the event $\{(w, U') \in \Omega | w, U' : w \in U\}$ that the value of $X$ lies in $U$ (and some $U'$ is observed). This latter subset of $\Omega$ is simply denoted by $U$, and is not strictly distinguished from $U$ as an event in the sample space $W$. Figure 1 illustrates these definitions for a three-element complete-data space $W = \{w_1, w_2, w_3\}$. The elements of $\Omega(W)$ correspond to the unfilled cells in this graphical representation. For $U = \{w_2, w_3\}$ the events $O_U$ and $U$ (as a subset of $\Omega$) are outlined.

A distribution $P$ on $\Omega$ is parameterized by the parameters $\theta$ defining the marginal distribution on $W$, and parameters
$$\lambda_{w,U} := P((w, U) | w), \qquad (w, U) \in \Omega,$$
defining the coarsening process.

EXAMPLE 2.2. Table 1 specifies distributions $P^{(i)}$, $i = 1, 2, 3$, on $\Omega(\{w_1, w_2, w_3\})$ through parameters $\theta^{(i)}$ on $W$ and conditional probabilities $\lambda^{(i)}$. For $w$ with $P_\theta(w) = 0$ parameters $\lambda_{w,U}$ are shown in brackets. Changing these parameters to arbitrary other values just leads to a different version of the conditional distribution of coarse observations given complete data, and has no influence on the joint distribution.

TABLE 1
*Parameters for distributions $P^{(1)}, P^{(2)}, P^{(3)}$*

| | | $\theta^{(i)}$ | $\lambda^{(i)}$ | | | | | | |
| --- | --- | --- | --- | --- | --- | --- | --- | --- | --- |
| | | | $\{w_1\}$ | $\{w_2\}$ | $\{w_3\}$ | $\{w_1,w_2\}$ | $\{w_1,w_3\}$ | $\{w_2,w_3\}$ | $\{w_1,w_2,w_3\}$ |
| $i=1$ | $w_1$ | 0 | [1/3] | | | [1/3] | [0] | | [1/3] |
| | $w_2$ | 1 | | 0 | | 1/3 | | 1/3 | 1/3 |
| | $w_3$ | 0 | | | [1/3] | | [0] | [1/3] | [1/3] |
| $i=2$ | $w_1$ | 1/2 | 0 | | | 2/3 | 0 | | 1/3 |
| | $w_2$ | 0 | | [1/3] | | [2/3] | | [0] | [0] |
| | $w_3$ | 1/2 | | | 0 | | 0 | 2/3 | 1/3 |
| $i=3$ | $w_1$ | 1/3 | 1/3 | | | 1/3 | 0 | | 1/3 |
| | $w_2$ | 1/3 | | 0 | | 1/3 | | 1/3 | 1/3 |
| | $w_3$ | 1/3 | | | 1/3 | | 0 | 1/3 | 1/3 |

As in this example, we generally assume that parameters $\lambda_{w,U}$ exist even when $P_\theta(w) = 0$ (rather than treating them as undefined), because in that way the parameter space $\Lambda^n$ for the $\lambda$-parameters does not depend on $\theta$:

$$\Lambda^n := \left\{ (\lambda_{w,U})_{w \in W, U \subseteq W : w \in U} | \lambda_{w,U} \in [0,1]; \forall w : \sum_{U : w \in U} \lambda_{w,U} = 1 \right\}.$$

Any subset $\Sigma \subseteq \Delta^n \times \Lambda^n$ is called a *coarse data model*. Such a model encodes assumptions both on the underlying complete data distribution and on the coarsening process. The complete-data model underlying $\Sigma$ is

$$\Theta = \{\theta \in \Delta^n | \exists \lambda : (\theta, \lambda) \in \Sigma\}.$$

We sometimes write $\Sigma(\Theta)$ for $\Sigma$ to emphasize the underlying complete-data model. We denote with $\Sigma_{sat}(\Theta) = \Theta \times \Lambda^n$ the saturated coarsening model with underlying $\Theta$.

A sample of coarse data items $\mathcal{U} = U_1, \ldots, U_N$ ($U_i \subseteq W$) is interpreted with respect to a coarse data model as observations of events $O_{U_i}$ in $\Omega$, and gives rise to the *observed-data likelihood* for $\theta$ and $\lambda$,

$$(1) \qquad L_{\text{OD}}(\theta, \lambda | \mathcal{U}) := \prod_{i=1}^{N} P_{\theta,\lambda}(O_{U_i}).$$

When ignoring the coarsening process, the data items $U_i$ are simply interpreted as subsets of $W$ and give rise to the *face-value likelihood* [4] for $\theta$,

$$(2) \qquad L_{\text{FV}}(\theta|\mathcal{U}) := \prod_{i=1}^{N} P_\theta(U_i).$$

**3. Ignorability.** The question of ignorability is under what conditions inferences about $\theta$ based on the face-value likelihood will be the same as obtained from the observed-data likelihood. These conditions will depend on the inference methods used [15]. Here we focus on the problem of ignorability for likelihood-based inference, with special emphasis on maximum likelihood estimation, which plays an important role in practice through the widespread use of the EM algorithm [6, 14].

For likelihood-based inference about $\theta$, the observed-data likelihood will typically be reduced to the *profile-likelihood*

$$(3) \qquad L_{P,\Sigma}(\theta|\mathcal{U}) := \max_{\lambda \,:\, (\theta,\lambda) \in \Sigma} L_{\text{OD}}(\theta, \lambda|\mathcal{U}).$$

To make the profile-likelihood well defined for all $\theta$, we restrict ourselves to models $\Sigma$ for which $\{\lambda|(\theta,\lambda) \in \Sigma\}$ is closed for every $\theta \in \Theta$, so that the maximum in (3) is attained. In our notation we make explicit that the profile-likelihood is not only a function of $\theta$ and $\mathcal{U}$, but also of the coarse data model $\Sigma$.

Moving from the observed-data likelihood to the profile-likelihood enables us to treat inference both with and without taking the coarsening process into account as inference with a likelihood function of only the parameter of interest, $\theta$. In particular, we obtain succinct formulations of ignorability questions: under what conditions on $\Sigma$ are likelihood ratios $L_{P,\Sigma}(\theta)/L_{P,\Sigma}(\theta')$ and $L_{\text{FV}}(\theta)/L_{\text{FV}}(\theta')$ equal for all $\theta, \theta'$; under what conditions are $L_{P,\Sigma}$ and $L_{\text{FV}}$ maximized by the same values $\theta \in \Theta$?

In the following we formulate the *car* and parameter distinctness assumptions as such modeling assumptions on $\Sigma$. In the case of *car* it turns out that we must distinguish two different versions.

DEFINITION 3.1. The data is *weakly coarsened at random* (*w-car*) according to $P_{\theta,\lambda}$, if for all $U \subseteq W$ and all $w, w' \in U$

$$(4) \qquad P_\theta(w) > 0, \ P_\theta(w') > 0 \implies \lambda_{w,U} = \lambda_{w',U}.$$

DEFINITION 3.2. The data is *strongly coarsened at random* (*s-car*) according to $P_{\theta,\lambda}$, if for all $U \subseteq W$ and all $w, w' \in U$

$$(5) \qquad \lambda_{w,U} = \lambda_{w',U}.$$

The difference between weak and strong *car*, thus, is that *s-car* also imposes a restriction on conditional probabilities $P_{\theta,\lambda}(O_U|w)$ when $P(w) = 0$. This is the version of *car* used by Gill, van der Laan and Robins [7] for categorical data. Underlying this version of *car* is the notion of *car* being a condition on the coarsening mechanism alone, which must be formulated without reference to the underlying complete-data distribution. Weak *car*, on the other hand, appears to be the more appropriate version when *car* is seen as a condition on the joint distribution of complete and coarsened data.

Gill, van der Laan and Robins ([7], page 274) also give a definition for *car* in general sample spaces. In contrast to their definitions in the discrete setup, that definition reduces for finite sample spaces to *w-car*, not *s-car*. They pose as an open problem whether (in the terminology established by our preceding definitions) it is always possible to turn a *w-car* model into an *s-car* model by a suitable setting of the $\lambda_{w,U}$-parameters for those $w$ with $P_\theta(w) = 0$. Our next example shows that this is not the case.

EXAMPLE 3.3. All distributions in Table 1 are *w-car*, but only $P^{(1)}$ and $P^{(3)}$ are *s-car*: to check the *w-car* condition it only is necessary to verify that all unbracketed $\lambda_{w,U}$ in a column are pairwise equal. For *s-car* also equality of the bracketed parameters is required. This condition is violated in the last two columns for $P^{(2)}$. Moreover, it is not possible to replace the bracketed $\lambda^{(2)}$-values with different conditional probabilities in a way that *s-car* is satisfied, because the conditional probabilities for the observations $\{w_1, w_2\}$, $\{w_2, w_3\}$ and $\{w_1, w_2, w_3\}$ would have to add up to 5/3.

In the following we write *car* when we wish to refer uniformly to both versions of *car*, for example, in definitions that can be analogously given for both versions, or in statements that hold for both versions.

When $P_{\theta,\lambda}$ satisfies *car* we denote parameters $\lambda_{w,U}$ simply with $\lambda_U$. In the case of *w-car* this denotes the parameter $\lambda_{w,U}$ common for all $w$ of positive probability. When $P_\theta(U) = 0$, then $\lambda_U$ is not well defined for *w-car*. We denote with $\Sigma_{car}(\Theta)$ the subset of $\Sigma_{sat}(\Theta)$ consisting of those parameters according to which the data is *car*. For $\theta \in \Theta$ we denote with $\Lambda_{w\text{-}car}(\theta)$ the set of $\lambda \in \Lambda^n$ that satisfy (4). Thus, $\Sigma_{w\text{-}car}(\Theta) = \{(\theta, \lambda) | \theta \in \Theta, \lambda \in \Lambda_{w\text{-}car}(\theta)\}$. From Definition 3.1 it follows that $\text{support}(P_\theta) \subseteq \text{support}(P_{\theta'})$ implies $\Lambda_{w\text{-}car}(\theta) \supseteq \Lambda_{w\text{-}car}(\theta')$. For *s-car* we can simply define the set $\Lambda_{s\text{-}car}$ of coarsening parameters that satisfy (5), and have $\Sigma_{s\text{-}car}(\Theta) = \Theta \times \Lambda_{s\text{-}car}$.

The following definition provides an important alternative characterization of *w-car*.

DEFINITION 3.4. $P_{\theta,\lambda}$ satisfies the *fair evidence condition* if for all $w, U$ with $w \in U$,

(6) $\qquad P_{\theta,\lambda}(O_U) > 0 \implies P_{\theta,\lambda}(w|O_U) = P_\theta(w|U).$

The fair evidence condition is necessary to justify updating a probability distribution by conditioning when an observation is made that establishes the actual state to be a member of $U$ [8]. We now obtain:

THEOREM 3.5. *The following are equivalent for $P_{\theta,\lambda}$:*

(a) *$P_{\theta,\lambda}$ satisfies w-car.*
(b) *$P_{\theta,\lambda}$ satisfies the fair evidence condition.*
(c) *For all $w, U$ with $w \in U$ and $P_\theta(w) > 0$,*

$$P_{\theta,\lambda}(O_U|w) = P_{\theta,\lambda}(O_U)/P_\theta(U).$$

PROOF. (a)$\Rightarrow$(b): If $P_{\theta,\lambda}(O_U|w) = P_{\theta,\lambda}(O_U|w')$ for all $w, w' \in U$ with $P_\theta(w), P_\theta(w') > 0$, then this value is equal to $P_{\theta,\lambda}(O_U|U)$. Assume that $P_\theta(w) > 0$ [otherwise there is nothing to show for (6)]. Using $P_{\theta,\lambda}(U|O_U) = 1$, then $P_{\theta,\lambda}(w|O_U) = P_{\theta,\lambda}(O_U|w)P_\theta(w)/P_{\theta,\lambda}(O_U) = P_{\theta,\lambda}(O_U|U)P_\theta(w)/P_{\theta,\lambda}(O_U) = P_{\theta,\lambda}(U|O_U)P_\theta(w)/P_\theta(U) = P_\theta(w|U)$.

(b)$\Rightarrow$(c): Let $w \in U$ with $P_\theta(w) > 0$. Then $P_{\theta,\lambda}(O_U|w) = P_{\theta,\lambda}(w|O_U) \times P_{\theta,\lambda}(O_U)/P_\theta(w) = P_\theta(w|U)P_{\theta,\lambda}(O_U)/P_\theta(w) = P_{\theta,\lambda}(O_U)/P_\theta(U)$.

(c)$\Rightarrow$(a): Obvious. $\square$

EXAMPLE 3.6. To check the fair evidence condition for the distributions of Table 1, one has to verify that for each observation $O_U$, normalizing all nonbracketed entries in the $\lambda$-column for $O_U$ yields the conditional distribution of $P_\theta$ on $U$.

One might suspect that one can also obtain a "strong fair evidence condition" by considering the normalization of both the bracketed and the unbracketed $\lambda$-entries, and that this strong version of the fair evidence condition would be equivalent to *s-car*. However, already for $P^{(1)}$ (which is *s-car*), we see that for $U = \{w_1, w_2\}$ the normalization of the column for $O_U$ gives $(1/2, 1/2)$ on $U$, which is not $P_\theta(\cdot|U)$.

Gill, van der Laan and Robins ([7], page 260) claim the equivalence of the fair evidence condition and *s-car*. However, as our results show, fair evidence is equivalent to *w-car*, not *s-car*. (The error in the proof of Gill, van der Laan and Robins [7] lies in an (implicit) application of Bayes rule to conditioning events of zero probability.) A correct proof of the equivalence (a)$\Leftrightarrow$(b) also is given by Grünwald and Halpern [8]. We consider the equivalence with the fair evidence condition to be an important point in favor of *w-car* as opposed to *s-car*.

Weak and strong *car* are modeling assumptions that identify certain coarse data distributions for inclusion in our model. The second condition usually required for ignorability, parameter distinctness, on the other hand, is a global condition on the structure of the coarse data model.

DEFINITION 3.7. A coarse data model $\Sigma$ satisfies *parameter distinctness* (*pd*) iff $\Sigma = \Theta \times \Lambda$ for some $\Theta \subseteq \Delta^n, \Lambda \subseteq \Lambda^n$.

From *car* and *pd* ignorability for likelihood-based inference can be derived. We next restate Rubin's proof of this result, in a way that clearly separates the contributions made by *car* and *pd*. To begin, assume that $\Sigma \subseteq \Sigma_{car}$, and let $(\theta, \lambda) \in \Sigma$. Let $\mathcal{U}$ be a sample with $P_\theta(U_i) > 0$ for $i = 1, \ldots, N$. Now

$$L_{\text{OD}}(\theta, \lambda | \mathcal{U}) = \prod_{i=1}^N P_{\theta,\lambda}(O_{U_i}) = \prod_{i=1}^N \sum_{w \in U_i} P_{\theta,\lambda}((w, U_i))$$

$$= \prod_{i=1}^N \lambda_{U_i} \sum_{w \in U_i} P_\theta(w) = \prod_{i=1}^N \lambda_{U_i} P_\theta(U_i).$$

Thus

(7) $$L_{P,\Sigma}(\theta | \mathcal{U}) = c_\Sigma(\theta, \mathcal{U}) L_{\text{FV}}(\theta | \mathcal{U}),$$

where

(8) $$c_\Sigma(\theta, \mathcal{U}) := \max_{\lambda : (\theta, \lambda) \in \Sigma} \prod_{i=1}^N \lambda_{U_i}.$$

Now assume, too, that *pd* holds, that is, $\Sigma = \Theta \times \Lambda$. The right-hand side of (8) then simply becomes $\max_{\lambda \in \Lambda} \prod_{i=1}^N \lambda_{U_i}$, which no longer depends on $\theta$. $L_{P,\Sigma}$ and $L_{\text{FV}}$, thus, differ only by a constant, so that inferences based on likelihood ratios of $L_{\text{FV}}$ are justified.

This derivation also provides the answer to a somewhat subtle question that arises out of our analysis so far: we have assumed throughout that the coarse data will be analyzed correctly in the coarse data space $\Omega$ using the observed-data likelihood $L_{\text{OD}}$. However, interpreting the data in $\Omega$ means that we still are dealing with coarse data, because it now is seen to consist of observations of subsets $O_U$ of $\Omega$, not of complete observations $(w, U) \in \Omega$. The question then is whether we have gained anything: $L_{\text{OD}}$ really is nothing but the face-value likelihood of incomplete data in the more sophisticated complete-data space $\Omega$. Do we thus have to build a second-order coarse data model on top of $\Omega$, and so on? The answer is no, because the coarsening process that turns complete data $(w, U)$ from $\Omega$ into coarse observations $O_U$ always is ignorable: in the second-order coarsening model we have $\lambda_{(w,U'),O_U} = 1$ iff $U' = U$, which means that here the data is *car*, and the factor $c(\theta, \mathcal{U})$ in (7) is always equal to 1.

How can this ignorability result be used in practice? In most cases it is appealed to simply by stating that the *car* and *pd* assumptions are made, and that this justifies the use of the face-value likelihood. This, however, is a rather incomplete justification, because *car* and *pd* together are not well-defined modeling assumptions that determine a unique coarse data model $\Sigma$. To make the *car* and *pd* assumptions only means to assume that the coarse data model $\Sigma$ is some subset of $\Sigma_{car}(\Theta)$, and has product form $\Theta' \times \Lambda'$. In the case of *w-car*, nontrivial further modeling assumptions may have to be made to ensure that *pd* holds, because $\Sigma_{w\text{-}car}(\Theta)$ itself usually is not a product. The following example illustrates the consequences for likelihood-based inferences under *w-car*. From now on we write $L_{P,car}$ and $c_{car}$ for $L_{P,\Sigma_{car}(\Theta)}$, respectively $c_{\Sigma_{car}(\Theta)}$, and similarly for $\Sigma_{sat}(\Theta)$. The underlying $\Theta$ will always be clear from the context.

EXAMPLE 3.8. Let $\theta^{(i)}$, $i = 1,2,3$, be as in Example 2.2. Let $\mathcal{U}$ be a sample consisting of $U_1 = \{w_1, w_2\}$, $U_2 = \{w_2, w_3\}$, $U_3 = \{w_1, w_2, w_3\}$. It is readily verified that for $i = 1, 2, 3$,

$$c_{w\text{-}car}(\theta^{(i)}, \mathcal{U}) = \lambda_{U_1}^{(i)} \cdot \lambda_{U_2}^{(i)} \cdot \lambda_{U_3}^{(i)},$$

that is, the coarsening parameters given in Table 1 maximize $\lambda_{U_1} \cdot \lambda_{U_2} \cdot \lambda_{U_3}$ over all parameters in $\Lambda_{w\text{-}car}(\theta^i)$. It also follows immediately that $c_{s\text{-}car}(\theta^{(i)}, \mathcal{U}) = (1/3)^3 = 1/27$. With these $c_{car}$-values one obtains the likelihood values shown in Table 2.

The first two columns of Table 2 show that likelihood ratios of $L_{FV}$ and $L_{P,w\text{-}car}$ do not coincide. Also the weaker ignorability condition of identical likelihood maxima does not apply: $L_{P,w\text{-}car}$ has the two maxima $P^{(1)}$ and $P^{(2)}$ (Theorem 4.4 below will show that these are indeed global maxima of $L_{P,w\text{-}car}$), but of these only $P^{(1)}$ also maximizes $L_{FV}$. It is not surprising that ignorability here cannot be established on the basis of *w-car* alone, because $\Sigma_{w\text{-}car}$ does not satisfy *pd*, and hence the factors $c_{w\text{-}car}(\theta, \mathcal{U})$ in (7) are different for different $\theta$. However, in Section 4.1 we will see that even on the basis of *w-car* alone a useful ignorability result can be obtained.

The *s-car* assumption, on the other hand, yields ignorability in the strong sense of equal likelihood ratios, because $\Sigma_{s\text{-}car} = \Theta \times \Lambda_{s\text{-}car}$ satisfies *pd*.

We thus obtain the following picture on the interdependence between the *car* and *pd* assumptions: *s-car* as the only modeling assumption on the coarsening process implies *pd*. To obtain ignorability, it therefore is sufficient to stipulate *s-car*. When one stipulates *w-car* as a modeling assumption, then additional assumptions are required to make the resulting model also satisfy *pd*. It must be realized that *pd* is itself not a well-defined modeling assumption, because it does not identify any particular subset of distributions for inclusion in the model. A joint assumption of *w-car* and *pd* only is possible if suitable further restrictions on either the complete-data model $\Theta$ or on the coarsening process are made. One possible restriction on $\Theta$ is to assume a fixed set of support for the distributions $P_\theta$. If, for example, $\Theta \subseteq \{\theta | \text{support}(P_\theta) = W\}$, then $\Sigma_{w\text{-}car}(\Theta)$ has *pd*. However, in most cases it is not possible to determine a priori the set of support of a categorical data distribution under investigation, and hence models allowing for different sets of support have to be used.

TABLE 2
*Likelihood values*

| $i$ | $L_{FV}(\theta^{(i)}|\mathcal{U})$ | $L_{P,w\text{-}car}(\theta^{(i)}|\mathcal{U})$ | $L_{P,s\text{-}car}(\theta^{(i)}|\mathcal{U})$ |
|---|---|---|---|
| 1 | 1 | $1 \cdot 1/27$ | $1 \cdot 1/27$ |
| 2 | 1/4 | $1/4 \cdot 4/27$ | $1/4 \cdot 1/27$ |
| 3 | 4/9 | $4/9 \cdot 1/27$ | $4/9 \cdot 1/27$ |

A further assumption one can make on the coarsening mechanism is that the data is *completely coarsened at random* (*ccar*) [9]. We do not give the precise definitions here, but only note that $\Sigma_{ccar}(\Theta) \subseteq \Sigma_{s\text{-}car}(\Theta)$ for any $\Theta$, and that $\Sigma_{ccar}(\Theta)$ has *pd*. Thus *ccar*, too, guarantees ignorability when it is the only modeling assumption on the coarsening mechanism. However, *ccar* is considered to be an unrealistically strong assumption for most applications.

**4. Ignorability without parameter distinctness.** In the preceding section we have seen that standard ignorability conditions cannot be established from the *w-car* assumption alone, because $\Sigma_{w\text{-}car}$ does not have *pd*. In this section we pursue the question whether some ignorability results can nevertheless be obtained from *w-car*. It turns out that in the case of the saturated complete-data model $\Theta = \Delta^n$ a fairly strong ignorability result for maximum likelihood inference can be obtained (Section 4.1). For nonsaturated complete-data models *s-car* is needed for ignorability. However, with nonsaturated models *car* becomes a testable assumption that, based on the observed data, may have to be rejected against the *not-car* alternative (Section 4.2).

The following simple lemma pertains to both saturated and nonsaturated complete-data models. For the formulation of the lemma we introduce the notation $c_{w\text{-}car}(V, \mathcal{U})$ for $c_{w\text{-}car}(\theta, \mathcal{U})$, where $\theta \in \Theta$ is such that $\text{support}(P_\theta) = V \subseteq W$. As $c_{w\text{-}car}(\theta, \mathcal{U})$ depends on $\theta$ only through $\text{support}(P_\theta)$, this is unambiguous. Identity (7) then becomes

$$(9) \qquad L_{P,w\text{-}car}(\theta | \mathcal{U}) = c_{w\text{-}car}(\text{support}(P_\theta), \mathcal{U}) L_{\text{FV}}(\theta | \mathcal{U}).$$

The following lemma is immediate from the definitions.

LEMMA 4.1. *Let* $V \subseteq V' \subseteq W$. *Then* $c_{w\text{-}car}(V, \mathcal{U}) \geq c_{w\text{-}car}(V', \mathcal{U})$.

4.1. *The saturated model.* For this section let $\Theta = \Delta^n$ be the saturated complete-data model. We immediately obtain a weak ignorability result.

THEOREM 4.2. *Let $\hat{\theta} \in \Delta^n$ be a local maximum of $L_{\mathrm{FV}}(\cdot|\mathcal{U})$. Then $\hat{\theta}$ is also a local maximum of $L_{P,w\text{-}car}(\cdot|\mathcal{U})$.*

PROOF. Let $\hat{\theta}$ be a local maximum of $L_{\mathrm{FV}}(\theta|\mathcal{U})$. There exists a neighborhood $\tilde{\Theta}$ of $\hat{\theta}$ such that for every $\tilde{\theta} \in \tilde{\Theta}$ we have $\mathrm{support}(P_{\tilde{\theta}}) \supseteq \mathrm{support}(P_{\hat{\theta}})$. With (9) and Lemma 4.1 the theorem then follows. $\square$

We next show that local maxima of $L_{\mathrm{FV}}$ are, in fact, global maxima of $L_{P,w\text{-}car}$, thus establishing ignorability in a strong sense for maximum likelihood inference in the saturated model. For the characterization of the maxima of $L_{P,w\text{-}car}$ the following definitions are needed. The terminology here is adopted from [5].

DEFINITION 4.3. Let $\Omega(W)$ be as in Definition 2.1. We denote with $\mathcal{O}$ the partition $\{O_U | \varnothing \neq U \subseteq W\}$ of $\Omega$. Let $m$ be a probability distribution on $\mathcal{O}$ and let $P_\theta$ be a probability distribution on $W$. We say that $m$ and $P_\theta$ are *compatible*, written $m \sim P_\theta$, if there exist parameters $\lambda \in \Lambda^n$, such that $P_{\theta,\lambda}(O_U) = m(O_U)$ for all $O_U \in \mathcal{O}$. We say that $m$ and $P_\theta$ are *car-compatible* (written $m \sim_{car} P_\theta$) if there exists such a $\lambda \in \Lambda_{car}(\theta)$.

THEOREM 4.4. *Let $\mathcal{U}$ be a set of data, and let $m$ be the empirical distribution induced by $\mathcal{U}$ on $\mathcal{O}$. For $\hat{\theta} \in \Delta^n$ with $\mathrm{support}(P_{\hat{\theta}}) = V \subseteq W$ the following are equivalent:*

(a) $m \sim_{w\text{-}car} P_{\hat{\theta}}$.
(b) $\hat{\theta}$ *is a global maximum of* $L_{P,w\text{-}car}(\theta|\mathcal{U})$ *in* $\Delta^n$.
(c) $L_{\mathrm{FV}}(\hat{\theta}|\mathcal{U}) > 0$, *and $\hat{\theta}$ is a local maximum of $L_{\mathrm{FV}}(\theta|\mathcal{U})$ within $\{\theta \in \Delta^n | \mathrm{support}(P_\theta) = V\}$.*

Theorem 4.4 establishes ignorability for maximum likelihood inference in a slightly different version from our original formulation in Section 3: it is not the case that $L_{P,w\text{-}car}$ and $L_{\mathrm{FV}}$ are (globally) maximized by the same $\theta \in \Theta$; however, maximization of $L_{\mathrm{FV}}$ will nevertheless produce the desired maxima of $L_{P,w\text{-}car}$, and, moreover, only a local maximum of $L_{\mathrm{FV}}$ must be found.

The proof of the theorem is preceded by two lemmas. The first one characterizes maxima of the observed-data likelihood in the saturated coarse data model.

LEMMA 4.5. *Let $\mathcal{U}$ and $m$ be as in Theorem 4.4. For $\hat{\theta} \in \Delta^n$ then the following are equivalent*:

(i) $m \sim P_{\hat{\theta}}$.
(ii) $\hat{\theta}$ is a global maximum of $L_{P,sat}(\theta|\mathcal{U})$.

PROOF. The likelihood $L_{\text{OD}}(\theta, \lambda|\mathcal{U})$ only depends on the marginal of $P_{\theta,\lambda}$ on $\mathcal{O}$, and is thus maximized whenever this marginal agrees with the empirical distribution.

The equivalence (i)⇔(ii) follows, because for every empirical distribution $m$ there exists at least one parameter $(\hat{\theta}, \hat{\lambda}) \in \Sigma_{sat}(\Delta^n)$ such that the marginal of $P_{\hat{\theta},\hat{\lambda}}$ on $\mathcal{O}$ is $m$. □

LEMMA 4.6. *The following are equivalent:*

(i) $m \sim_{w\text{-}car} P_\theta$.
(ii) *For all* $w \in W : P_\theta(w) > 0 \Rightarrow \sum_{U\,:\,w\in U} \frac{m(O_U)}{P_\theta(U)} = 1$.

The proof follows easily from Theorem 3.5.

PROOF OF THEOREM 4.4. (a)⇒(b): $m \sim_{w\text{-}car} P_{\hat{\theta}}$ trivially implies $m \sim P_{\hat{\theta}}$. By Lemma 4.5 $L_{P,sat}$ is maximized by $\hat{\theta}$. Also, $L_{P,sat}(\theta|\mathcal{U}) \geq L_{P,w\text{-}car}(\theta|\mathcal{U})$ with equality for $\theta = \hat{\theta}$. Hence $\hat{\theta}$ maximizes $L_{P,w\text{-}car}$.

(b)⇒(c): Immediate from (9).

(c)⇒(a): Recall that $\theta \in \Delta^n$ is written as $\theta = (p_1, \ldots, p_n)$ with $p_i = P_\theta(w_i)$. Let $D := \{\theta \in \Delta^n | L_{\text{FV}}(\theta|\mathcal{U}) > 0\}$. For $\theta \in D$ then

$$\frac{1}{N}\log L_{\text{FV}}(\theta|\mathcal{U}) = \sum_{U \subseteq W\,:\,m(O_U)>0} m(O_U) \log P_\theta(U)$$

$$= \sum_{U \subseteq W\,:\,m(O_U)>0} m(O_U) \log\left(\sum_{i\,:\,w_i\in U} p_i\right).$$

This is differentiable on $D$, and with $\mathcal{U}(w_i) := \{U \subseteq W | m(O_U) > 0, w_i \in U\}$,

$$\frac{\partial (1/N) \log L_{\text{FV}}(\theta|\mathcal{U})}{\partial p_i} = \sum_{U \in \mathcal{U}(w_i)} m(O_U)\left(\sum_{j\,:\,w_j\in U} p_j\right)^{-1} = \sum_{U \in \mathcal{U}(w_i)} \frac{m(O_U)}{P_\theta(U)}$$

[the sum on the right-hand side being interpreted as 0 when $\mathcal{U}(w_i) = \varnothing$]. For $\hat{\theta}$ as in (c) we have that $S(V) := \{\theta \in \Delta^n | \text{support}(P_\theta) = V\} \subseteq D$, and the gradient of $(1/N) \log L_{\text{FV}}(\theta|\mathcal{U})$ is orthogonal to $S(V)$ at $\hat{\theta}$. This can be expressed as the condition that for every $\theta' = (p'_1, \ldots, p'_n) \in S(V)$

$$\sum_{i=1}^n \left(\sum_{U \in \mathcal{U}(w_i)} \frac{m(O_U)}{P_{\hat{\theta}}(U)}\right)(\hat{p}_i - p'_i) = 0,$$

which is equivalent to

$$\sum_{i:\,w_i\in V}\left(\sum_{U\in\mathcal{U}(w_i)}\frac{m(O_U)}{P_{\hat\theta}(U)}\right)\hat p_i = \sum_{i:\,w_i\in V}\left(\sum_{U\in\mathcal{U}(w_i)}\frac{m(O_U)}{P_{\hat\theta}(U)}\right)p_i'.$$

This implies that $\sum_{U\in\mathcal{U}(w_i)} m(O_U)/P_{\hat\theta}(U)$ is a constant $k$ that does not depend on $w_i$, and furthermore

$$k = k\sum_{i:\,w_i\in V}\hat p_i = \sum_{i:\,w_i\in V}\sum_{U\in\mathcal{U}(w_i)}\frac{m(O_U)}{P_{\hat\theta}(U)}\hat p_i = \sum_{U:\,m(O_U)>0}\sum_{i:\,w_i\in U}\frac{m(O_U)}{P_{\hat\theta}(U)}\hat p_i$$
$$= \sum_{U:\,m(O_U)>0} m(O_U) = 1.$$

Now (a) follows from Lemma 4.6. □

We remark that (a)–(c) in Theorem 4.4 also are equivalent to:

(d) $\hat\theta$ is a stationary point for the EM algorithm.

We do not give a formal proof here, but emphasize that in the current context we assume the saturated complete-data model, and thus in (d) also assume that the EM algorithm operates on the unrestricted parameter space $\Delta^n$. Then the equivalence (a)⇔(d) easily follows from the equivalence of *w-car* and the fair evidence condition.

For *s-car* one obtains the following analogue of Theorem 4.4.

THEOREM 4.7. *Let $\mathcal{U}, m$ be as in Theorem 4.4. For $\hat\theta\in\Delta^n$ the following are equivalent:*

(a) $m \sim_{s\text{-}car} P_{\hat\theta}$.
(b) $\hat\theta$ *is a global maximum of* $L_{P,s\text{-}car}(\theta|\mathcal{U})$ *in* $\Delta^n$.
(c) $\hat\theta$ *is a global maximum of* $L_{\mathrm{FV}}(\theta|\mathcal{U})$ *in* $\Delta^n$.

The equivalence (b)⇔(c) here is immediate from the equality of likelihood ratios of $L_{P,s\text{-}car}$ and $L_{\mathrm{FV}}$. The nontrivial implication is (c)⇒(a). It has (implicitly) been shown by Gill, van der Laan and Robins [7] in the proof of their first theorem.

EXAMPLE 4.8. Let $\mathcal{U}$ be as in Example 3.8. Then $m(O_{U_i}) = 1/3$ for $i=1,2,3$. $P^{(1)}$ and $P^{(2)}$ are *w-car* distributions with marginal $m$ on $\mathcal{O}$. By Theorem 4.4, $\theta^{(1)}$ and $\theta^{(2)}$ are global maxima of $L_{P,w\text{-}car}$. $P^{(1)}$ also is *s-car*, and therefore $\theta^{(1)}$ is a global maximum of $L_{P,s\text{-}car}$.

In the preceding example we found a single maximum of $L_{P,s\text{-}car}$, and two distinct maxima of $L_{P,w\text{-}car}$. Gill, van der Laan and Robins [7] showed that for every $m$ there exists $\theta$ with $m \sim_{s\text{-}car} P_\theta$, and $\theta$ is essentially unique [for any $\theta'$ with $m \sim_{s\text{-}car} P_{\theta'}$ it must be the case that $P_{\theta'}(U) = P_\theta(U)$ for all $U \in \mathcal{U}$]. Thus $L_{P,s\text{-}car}$ has an essentially unique maximum. For $w\text{-}car$ we obtain the following result on the existence of maxima of $L_{P,w\text{-}car}$.

THEOREM 4.9. *Let $\mathcal{U}$ and $m$ be as in Theorem 4.4. Let $V \subseteq W$ such that for all $U \subseteq W$*

(10) $$m(O_U) > 0 \implies V \cap U \neq \varnothing.$$

*Then there exists $\hat{\theta}$ with $\text{support}(P_{\hat{\theta}}) \subseteq V$ and $m \sim_{w\text{-}car} P_{\hat{\theta}}$.*

PROOF. From (10) it follows that $L_{\text{FV}}(\theta|\mathcal{U}) > 0$ for $\theta$ with $\text{support}(P_\theta) = V$. In particular, $L_{\text{FV}}(\theta|\mathcal{U})$ is not identically zero on the compact set $\{\theta \,|\, \text{support}(P_\theta) \subseteq V\}$ and attains a positive maximum at some $\hat{\theta}$. The theorem now follows from Theorem 4.4. $\square$

Theorem 4.4 in conjunction with Lemma 4.5 provides yet another ignorability result: maximization of $L_{\text{FV}}$ will yield a parameter $\hat{\theta}$ with $m \sim P_{\hat{\theta}}$, and thus a global maximum of $L_{P,sat}$. Thus, the use of the face-value likelihood instead of the observed-data likelihood also is justified when we assume the saturated coarse data model $\Sigma_{sat}(\Delta^n)$. In other words, ignorability holds when the coarsening process is treated as completely unknown (and the saturated model also is assumed for the underlying complete data). However, it turns out that ignorability is not really the issue here, as maximum likelihood solutions for the observed-data likelihood in the model $\Sigma_{sat}(\Delta^n)$ can be found directly without optimizing $L_{\text{FV}}$: Dempster [5] gives an explicit construction of the set $\{\theta \in \Delta^n \,|\, m \sim P_\theta\}$, which briefly is as follows.

Consider any ordering $w_{i_1}, w_{i_2}, \ldots, w_{i_n}$ of the elements of $W$. Now transform the coarse data $U_1, \ldots, U_N$ into a sample of complete-data items by interpreting $U_j$ as an observation of the first $w_{i_h}$ in the given ordering that is an element of $U_j$. Let $P_\theta$ be the empirical distribution of this completed sample. By considering all possible orderings of $W$, one obtains in this way distributions $P_{\theta_1}, \ldots, P_{\theta_{n!}}$ on $W$. The set $\{\theta \in \Delta^n \,|\, m \sim P_\theta\}$ now is the convex hull of all these $P_{\theta_i}$. Moreover, the empirical distribution of any completion of the data lies in the convex hull of the $P_{\theta_i}$.

It thus is very easy to directly determine some maximal likelihood solutions of $L_{P,sat}$, simply as the empirical distribution of an arbitrary completion of the data. An explicit representation of all solutions is obtained by computing all $P_{\theta_i}$.

The problem here is that the set $\{\theta \in \Delta^n | m \sim P_\theta\}$ typically will be very large (much larger than the set $\{\theta \in \Delta^n | m \sim_{car} P_\theta\}$), and therefore inferences based on the coarse data model $\Sigma_{sat}(\Delta^n)$ will be too weak for practical purposes. We thus see that making the *car* assumption really serves a second purpose besides justifying the use of the face-value likelihood: we need to make some assumptions on the coarsening mechanism, because otherwise our model will be too weak to support practically useful inferences.

Figure 2 summarizes some of our results in terms of our running example. Shown is the polytope $\Delta^3$ with the potential lines of the face-value likelihood $L_{\mathrm{FV}}(\cdot|\mathcal{U})$ for $\mathcal{U}$ as in Example 3.8. The two distributions $P^{(1)}, P^{(2)}$ are marked by circles. They correspond to nonzero maxima of $L_{\mathrm{FV}}$ relative to distributions with the same set of support. Marked as diamonds are the distributions obtained from the extremal data completions for the five possible orderings of $W$. Their convex hull is the set of $\theta$ compatible with $m$.

From the results of this section we can also retrieve Gill, van der Laan and Robins' [7] result that "*car* is everything," that is, the *car* assumption cannot be rejected against the *not-car* alternative based on any observed coarse data (assuming an underlying saturated complete-data model). This is because by Theorem 4.9 [and the corresponding result in [7] for *s-car*] there exists for any observed coarse data a *car* model with the observed marginal on $\mathcal{O}$. Gill, van der Laan and Robins [7] show that the same need not hold for infinite sample spaces. Further results on the nontestability of the *car* assumption in general sample spaces have been obtained by Cator [1]. In the next section we will see that *car* also becomes testable for finite sample spaces with a parametric complete-data model.

4.2. *Nonsaturated models.* Most of the preceding ignorability results are no longer valid when the complete-data model is not $\Delta^n$. Only the weak ignorability result of Theorem 4.2 can be retained for a wide class of complete-data models.

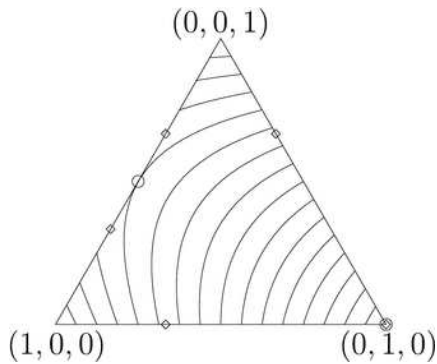

Fig. 2. *Summary of running example.*

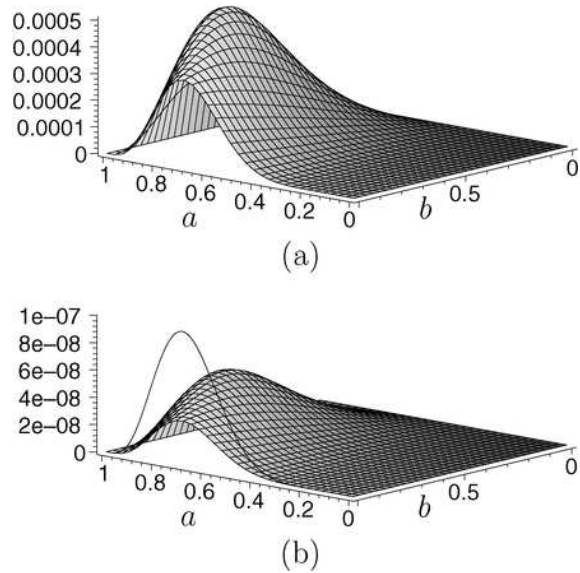

Fig. 3. $L_{\text{FV}}$ and $L_{P,w\text{-}car}$ in Example 4.12.

DEFINITION 4.10. A complete-data model $\{P_\theta | \theta \in \Theta\}$ is *support-continuous* if for all $\theta \in \Theta$ there exists a neighborhood $G_\theta \subseteq \Theta$ such that $\text{support}(P_{\theta'}) \supseteq \text{support}(P_\theta)$ for all $\theta' \in G_\theta$.

Virtually all natural parametric models are support-continuous. The proof of Theorem 4.2 actually has established the following:

THEOREM 4.11. *Let $\{P_\theta | \theta \in \Theta\}$ be a support-continuous complete-data model. Every local maximum $\hat{\theta} \in \Theta$ of $L_{\text{FV}}(\cdot|\mathcal{U})$ then is a local maximum of $L_{P,w\text{-}car}(\cdot|\mathcal{U})$.*

The following example shows that other results of Section 4.1 cannot be extended to parametric models.

EXAMPLE 4.12. Let $A$ and $B$ be two binary random variables. Let $W = \{AB, A\bar{B}, \bar{A}B, \bar{A}\bar{B}\}$, where, for example, $A\bar{B}$ represents the state where $A = 1$ and $B = 0$. We represent a probability distribution $P$ on $W$ as the tuple $(P(AB), P(A\bar{B}), P(\bar{A}B), P(\bar{A}\bar{B}))$. Define

$$\Theta = \{\theta = (a,b) | a \in [0,1], b \in [0,1]\},$$
$$P_\theta = (ab, a(1-b), (1-a)(1-b), (1-a)b).$$

Now assume that the data $\mathcal{U}$ consists of six observations of $A$ (i.e., the set $\{AB, A\bar{B}\}$), three observations of $B$, three observations of $\bar{B}$ and one observation of $\bar{A}\bar{B}$.

Figure 3(a) shows a plot of $L_{\mathrm{FV}}(\theta|\mathcal{U})$. We can numerically determine the unique maximum as $\hat\theta \approx (0.845, 0.636)$, which corresponds to $P_{\hat\theta} \approx (0.54, 0.31, 0.05, 0.1)$. Restricted to the subset $\Theta_1 := \{\theta \in \Theta | 0 < a < 1, b = 1\}$ a local maximum is attained at $\theta^* \approx (0.69, 1)$, which corresponds to $P_{\theta^*} \approx (0.69, 0, 0, 0.31)$.

The set $\Theta_1$ corresponds to the set of support $V = \{AB, \bar{A}\bar{B}\}$ of $P_\theta$, that is, $\theta \in \Theta_1 \Leftrightarrow \mathrm{support}(P_\theta) = V$. Similarly, the set $\Theta_2 := \{\theta \in \Theta | 0 < a < 1, 0 < b < 1\}$ contains the parameters $\theta$ that define distributions with full set of support $W$. $L_{P,w\text{-}car}(\theta|\mathcal{U})$, therefore, is given by multiplying $L_{\mathrm{FV}}(\theta|\mathcal{U})$ by $c_{w\text{-}car}(V,\mathcal{U})$ when $\theta \in \Theta_1$, and by $c_{w\text{-}car}(W,\mathcal{U})$ when $\theta \in \Theta_2$. For all $\theta \notin \Theta_1 \cup \Theta_2$ we obtain $L_{\mathrm{FV}}(\theta|\mathcal{U}) = 0$, so that further constants $c_{w\text{-}car}(V'|\mathcal{U})$ do not matter. The approximate values for the relevant constants are $c_{w\text{-}car}(V,\mathcal{U}) \approx 0.0003$ and $c_{w\text{-}car}(W,\mathcal{U}) \approx 0.0001$.

A plot of $L_{P,w\text{-}car}(\theta|\mathcal{U})$ as given by (9) is shown in Figure 3(b). Note the discontinuity at the boundary between $\Theta_1$ and $\Theta_2$ due to the different factors $c_{w\text{-}car}(V,\mathcal{U})$ and $c_{w\text{-}car}(W,\mathcal{U})$. It turns out that the global maximum now is $\theta^*$, rather than $\hat\theta$.

Theorem 4.4 allows us to analyze the situation more clearly. It is easy to see that $P = (9/13, 0, 0, 4/13)$ is a distribution that is $w\text{-}car$-compatible with the empirical distribution $m$ induced by $\mathcal{U}$. We find that $P = P_{\theta^*}$ for $\theta^* = (9/13, 1) = (0.6923, 1)$, which thus turns out to be the precise value of $\theta^*$ which initially was determined numerically. From Theorem 4.4 it now follows that $P_{\theta^*}$ has maximal $L_{P,w\text{-}car}$-likelihood score even within the class of all distributions on $W$, so that not only is $\theta^*$ a global maximum in $\Theta$, but no better solution can be found by changing the parametric complete-data model.

Under the $s\text{-}car$ assumption the maximum likelihood estimate is $\hat\theta$. Thus, the two versions of $car$ here lead to quite different inferences. There also is a fundamental difference with respect to testability: while $\theta^*$ is $w\text{-}car$-compatible with $m$, $\hat\theta$ is not $s\text{-}car$-compatible with $m$. Consequently, the $s\text{-}car$ hypothesis, but not the $w\text{-}car$ hypothesis, can be rejected against the unrestricted alternative $\Sigma_{sat}$ by a likelihood ratio test (when $m$ is induced by a sufficiently large sample).

We can summarize the results for nonsaturated models as follows: since $\Sigma_{s\text{-}car}(\Theta)$ satisfies $pd$ for any parametric model $\Theta$, ignorability for likelihood-based inference is guaranteed by $s\text{-}car$.

For $w\text{-}car$, even the weak ignorability condition that maximization of $L_{\mathrm{FV}}$ will give a maximum of $L_{P,w\text{-}car}$ does not hold. The apparent advantage of $s\text{-}car$ has to be interpreted with caution, however: whenever a maximum $\hat\theta$ of $L_{\mathrm{FV}}$ maximizes $L_{P,s\text{-}car}$, but not $L_{P,w\text{-}car}$, then $P_{\hat\theta}$ cannot be $s\text{-}car$-compatible with $m$. Loosely speaking, this means that we obtain ignorability for maximum likelihood inference through $s\text{-}car$ but not through $w\text{-}car$ only when the data contradicts the *s-car* assumption. The same data, on the other hand, might be consistent with *w-car*, but for inference under the *w-car* assumption the face-value likelihood has to be corrected with the $c_{w\text{-}car}$ factors.

**5. Conclusion.** We can summarize the results of Sections 3 and 4 as follows: ignorability for maximum likelihood inference and categorical data holds under any of the following four modeling assumptions: 1. The *w-car* assumption for the coarsening mechanism and additional assumptions, such that the resulting coarse data model satisfies *pd*. 2. The *s-car* assumption as the sole assumption on the coarsening process. 3. The saturated complete-data model and *w-car* as the sole assumption on the coarsening process. 4. The saturated model for both the complete data and the coarsening mechanism (but here there are more efficient ways of finding likelihood maxima than by maximizing the face-value likelihood). In particular, one must be aware of the fact that the joint assumption $car + pd$ is ambiguous and can be inconsistent. This is because *pd* is not a well-defined modeling assumption one is free to make, but a *model property* one has to ensure by other assumptions.

Overall the ignorability results obtained from *s-car* are somewhat stronger than those obtained from *w-car*. Points in favor of *w-car*, on the other hand, are its equivalence with the fair evidence condition, and the fact that it is invariant for different versions of the conditional distribution of observed (coarse) data. Furthermore, the *w-car* assumption can be consistent with a given parametric model and observed data when *s-car* is not (but not vice versa).

**Acknowledgments.** The author thanks James Robins and Richard Gill for valuable discussions that clarified the intricacies of the weak-*car*, strong-*car* relationship. The original motivation for this work was in part provided by Ian Pratt by suggesting a "fair evidence principle" for conditioning.

Institut for Datalogi
Aalborg Universitet
Fredrik Bajers Vej 7E
DK-9220 Aalborg Ø
Denmark
e-mail: jaeger@cs.auc.dk